%% file: main.tex
\documentclass[fleqn,10pt]{wlscirep}
\usepackage[utf8]{inputenc}
\usepackage[T1]{fontenc}

\usepackage{graphicx}
\usepackage{subcaption}
\usepackage[utf8]{inputenc}
\usepackage{dcolumn}
\usepackage{bm}
\usepackage{color}
\usepackage{amsmath} 
\usepackage{physics}
\usepackage{braket}
\usepackage{amssymb}
\usepackage{amsthm}
\usepackage{mathtools} 
\usepackage{amsfonts}
\usepackage{comment}
\usepackage{xcolor}
\usepackage{url}
\usepackage{algorithm}
\usepackage{algorithmic}
\usepackage{soul}
\usepackage{caption}
\usepackage{subcaption}
\usepackage{silence}

\theoremstyle{definition}
\newtheorem{definition}{Definition}
\theoremstyle{definition}

\newtheorem{lemma}{Lemma}

\newtheorem{theorem}{Theorem}
\newtheorem{proposition}{Proposition}
\newtheorem{corollary}{Corollary}

\DeclareMathOperator*{\argmax}{arg\,max}


\title{A Compressive Sensing Inspired Monte-Carlo Method for Combinatorial Optimization}

\author[1,*]{Baptiste Chevalier}
\author[1]{Shimpei Yamaguchi}
\author[1]{Wojciech Roga}
\author[1,2]{Masahiro Takeoka}
\affil[1]{Department of Electronics and Electrical Engineering, Keio University, Yokohama 223-8522, Japan}
\affil[2]{Advanced ICT Research Institute, National Institute of Information and Communications Technology (NICT), Koganei, Tokyo 184-8795, Japan}

\affil[*]{chevalier.baptiste@keio.jp}

\begin{abstract}
	In this paper, we present a Monte-Carlo Compressive Optimization algorithm, a new method for black-box combinatorial optimization, which only assumes that the objective is compressible. The method relies on random queries to the objective function in order to estimate generalized moments. Next, a greedy algorithm from compressive sensing is repurposed to find the global optimum when not overfitting to the samples. We provide numerical results giving evidence that our method is competitive by comparing it with digital annealing. Moreover, we give theoretical justification for the success of the algorithm and analyze its properties. The practicality of our algorithm is enhanced by the ability to tune the heuristic parameters to the available computational resources. An end-to-end open-source implementation is available to use our method.
	
\end{abstract}
\begin{document}
	
	\flushbottom
	\maketitle

	\thispagestyle{empty}

	\section*{Introduction}
	
	Optimization problems lie at the very core of many important modern-world applications, including engineering, finance, and physics. However, finding efficient ways to solve these problems remains an important concern, since computational time usually increases quickly with the size of the instance. For example, the minimization of the cost function in neural networks over the space of network parameters becomes harder as the network depth increases. Optimization is generally difficult, but in some cases one can overcome the difficulty by exploiting the structure of the cost function. If the cost function is continuous, convex, or differentiable, powerful methods are known, such as gradient-based methods including Stochastic Gradient Descent \cite{robbins_stochastic_1951} or gradient-free methods including COBYLA or Nelder-Mead \cite{nelder_simplex_1965,powell_direct_1994}. 

	However, combinatorial optimization, which we focus on in this paper, raises different difficulties. The solution space of combinatorial problems is the set of all bit-strings of a given length $N$, $\{0,1\}^N$, so that naive methods -- brute force, linear search -- have a cost exponential in $N$, and the difficulty of a given problem depends on the structure of its objective. When the objective has a recursive structure, for instance interactions between neighboring bits only, dynamic programming methods \cite{bellman_theory_1954,schuch_matrix_2010} solve it exactly, at a cost that grows with the range of the interactions, since longer-range interactions increase the dimension of the state that must be kept. Many challenging combinatorial tasks can be expressed as generalized Ising Hamiltonians \cite{lucas_ising_2014,kirkpatrick_optimization_1983,kadowaki_quantum_1998,ohzeki_quantum_2010}, whose ground state is NP-hard to find for general couplings \cite{barahona_computational_1982}; these tasks include NP-hard problems such as the knapsack problem or the traveling salesman problem. Moreover, in many cases the objective has no explicit form and can only be evaluated. In this black-box scenario, methods that rely on the structure of the objective, such as dynamic programming, cannot be applied, and the best-known methods are meta-heuristics, such as simulated annealing \cite{kirkpatrick_optimization_1983,delahaye_simulated_2019} or genetic algorithms.
	
	In this work, we propose a novel optimization algorithm for black-box combinatorial optimization tasks. The method proposed in this paper is based on the principles of compressive sensing \cite{candes_robust_2006,donoho_compressed_2006,baraniuk_simple_2008,foucart_basic_2013}, statistical compressive learning \cite{gribonval_compressive_2021}, and Monte-Carlo methods \cite{metropolis_monte_1949,metropolis_equation_1953}. Compressive sensing is a technique that allows for recovering a high-dimensional signal from a small number of linear measurements, given that the signal has a sparse representation in a known basis. The recovery procedure explicitly finds the sparse representation from which the high-dimensional signal is derived by a basis transformation. The number of measurements needed to recover the total signal scales as a polynomial function of the number of non-zero entries of the sparse representation, which is much smaller than the dimensionality of the signal. 
	Statistical compressive learning is a relatively new technique for extracting the core features of a distribution from a sample. The sample is mixed by a sketch function, and a recovery algorithm applied to the sketching image returns a function that matches the original distribution better than the estimates based on the sample directly. 
	Finally, Monte Carlo is a well-known approach that uses stochastic algorithms to determine properties of high-dimensional functions that are intractable with rigorous analytical approaches. The combinatorial optimization method we propose here, Monte Carlo Compressive Optimization (MCCO), combines attractive features of the three techniques to find the maximum of a function based on a sample drawn from its values. MCCO derives the properties of the maximum value of the function from sketches applied to random samples, even though the maximum may never be met during the sampling process. In fact, the method can be seen as approximating a model whose parameters are learned from the samples and which is then used in the optimization procedure. We provide evidence that the location of the maximum can survive substantial modeling error, and that the choice of the model can be loose. This approach can be useful for a wide range of applications, including reinforcement learning and graph problem optimization.
	
	The scope of this paper is black-box combinatorial optimization: the solver can only query the objective, and the only assumption made on it is that it is compressible. Although the exact class of functions this method handles best remains unknown, we argue that this scope is large. Here, compressibility means that there exists a short formula determining, or well-approximating, all values of the cost function, or one can find an embedding from the compressed description space to the optimization space that maps all the relevant information into the values of the cost function. We show numerical evidence that MCCO is applicable to compressible functions, although the compressed representation does not need to be known explicitly. The numerical analysis shows that MCCO, when properly tuned and applied to the compressible combinatorial problems we describe, is able to carry out the optimization and is competitive with more usual methods such as digital annealing \cite{aramon_physics-inspired_2019}.
	
	Additionally, we developed an open-source Python library named TrOMA \cite{baptiste_baptistechevtroma_2026} (Tractable Optimization Modeling Approach) that integrates the MCCO method described in this work. Most importantly, the library is end-to-end and the API is problem- and hardware-agnostic. This means that anyone can test our method easily on any combinatorial problem; the only requirement is to provide a cost function. The library supports several optimization backends, including classical dynamic programming and heuristic methods for quantum computers (compatible with real IBM quantum computers).
	
	The paper is structured as follows. The Methods introduce the MCCO algorithm for combinatorial optimization, followed by the description of the class of compressible combinatorial problems. The Numerical Results study this class through numerical simulations; our results include a comparison between MCCO and digital annealing. The Theoretical Analysis strengthens our results by justifying the algorithm theoretically in general terms. The Discussion addresses applications of the MCCO method to Reinforcement Learning and to the Ising problem. We also discuss the role quantum computers can play in optimization. In the Supplementary Information, we give supplementary numerical results (Section 1), discuss the coherence of the structured sketches (Section 2), verify Theorem 1 numerically (Section 3) and give the proofs on the effect of thresholding (Section 4). 
	
	\section*{Methods}
	
	\subsection*{The Monte-Carlo Compressive Optimization algorithm}
	\label{sec:algorithm}
	
	An optimization problem consists of finding the optimum (maximum or minimum) of a given real-valued function $f$, called the objective, cost, loss, or reward function. Consider $f : \{0,1\}^N \rightarrow \mathbb{R}^+$, a cost function over the set of all possible bit-strings of length $N$. After normalization, $f$ is the density of a discrete probability distribution, which we denote $\pi$. Our goal is to find the best input bit-string $x^*$ that maximizes the function $f$. We define it as\\
	
	\noindent{\bf Problem I}$\qquad\qquad x^*_I=\argmax f(x)$.\\
	
\begin{figure}[ht]
		\centering
		\includegraphics[width=1\linewidth]{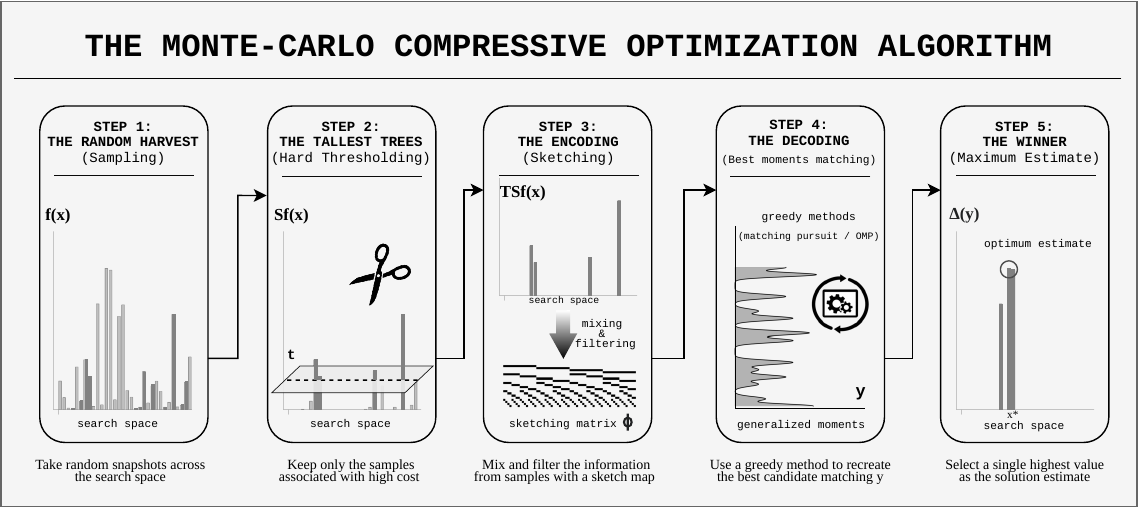}
		\caption{Schematic representation of Algorithm 1: the Monte Carlo Compressive Optimization. The algorithm proceeds through five stages: (1) Sampling of the search space; (2) Hard thresholding to isolate high-value samples; (3) Sketching, where the remaining samples are compressed into a compact sketch vector $y$; (4) Decoding via greedy reconstruction (e.g., matching pursuit); and (5) Estimation of the global maximum from the reconstructed vector.}
		\label{fig:scheme}
	\end{figure}

	\noindent Let us introduce a new algorithm for combinatorial optimization. The full scheme can be seen in Figure \ref{fig:scheme}.\\
	\noindent{\bf Algorithm 1} [Monte Carlo Compressive Optimization]\\
	\noindent \underline {Inputs:} sample size $n$, threshold parameter $t$, sketch function $\Phi$, decoding procedure $\operatorname{Decode}$.
	\begin{enumerate}
		\label{alg1}
		\item Compute a sample of $f$ called $Sf$. To do so, first draw a sample $\mathcal{I}$ of $n$ indices from the uniform distribution. Then apply the sampling operator $S_\mathcal{I}$, defined by $(S_\mathcal{I} f)(x) = f(x)$ if $x \in \mathcal{I}$ and $0$ otherwise.
		\item Apply a {\it hard thresholding} operator $T_t$ such that in $T_tSf$ all values $Sf(x)$ lower than a given $t$ become zero.$^*$
		\item Apply a {\it sketching map} (also known as a measurement map) $\Phi$ to get $y = \Phi(T_tSf)$. We call $y$ a {\it sketch vector}. It contains an empirical estimate of some {\it generalized moments} of the distribution $\pi$. $^{**}$
		\item Apply a decoding procedure $\operatorname{Decode}$ to $y$ and get $\tilde{f} = \operatorname{Decode}(y)$. The decoding procedure is taken from compressive sensing \textbf{greedy methods} (matching pursuit, orthogonal matching pursuit, \ldots).$^{***}$
		\item The largest 
		value of $\tilde{f}$ is used as the maximum estimate.  
	\end{enumerate}
	
	$^{*}$Throughout the analysis, the threshold is assumed to lie in the range $0 \le t \le f(x_2)$, 
	where $x_2 = \argmax_{x \neq x^*} f(x)$ is the second maximum of $f$, so that the gap of $f$ is left unchanged by $T_t$. 
	
	$^{**}$The Monte-Carlo method \cite{metropolis_monte_1949} can be used to approximate generalized moments similar to those of Ref.~\cite{gribonval_compressive_2021}. This is a well-known result from probabilistic computation theory (or ergodic theory) about a large number $M$ of samples drawn independently from a probability distribution $\pi$, 
	\begin{equation}\lim_{M\rightarrow\infty}
		\frac{1}{M}\sum_{x_{(i)}\sim \pi}g(x_{(i)}) =\sum_{x_{(i)}}g(x_{(i)})\pi(x_{(i)}) \approx \frac{2^N}{M}\sum_{x_{(i)}\sim\mathcal{U}}g(x_{(i)})\pi(x_{(i)}),
	\end{equation}
	where $2^N$ is the size of the domain.
	This means that a moment of a function $g$ of a variable $x$ with respect to a measure $\pi(x)$ can be approximated by sampling from the uniform distribution $\mathcal{U}$ and computing the value of the density $\pi(x)$, i.e. of $f(x)$ up to normalization, in those points.
	
	Regarding the choice of the sketch function, we know from compressive sensing theory that functions $\Phi$ that satisfy some properties such as the Restricted Isometry Property (RIP) are good candidates for sketch functions which should preserve important properties of the tested function -- in practice, often randomized linear maps are used. Structured deterministic measurement functions can also work \cite{foucart_basic_2013}. 
	For instance, in our past works \cite{jacob_franck-condon_2020,chevalier_compressed_2024}, structured measurement functions with different benefits were used. Moreover, the sketch function does not necessarily need to be linear; generalized moments of the distribution can be of different degrees. Introducing nonlinearities such as a hard thresholding process \cite{foucart_basic_2013} in practice gives us better performance than a random linear map alone. In addition, the choice of the sketch function directly impacts the optimization step within the recovery algorithm. The freedom in the choice of the sketch function allows for flexibility in the algorithm complexity, which can be tailored to match one's available computational resources. Indeed, this turned out to be an essential element of practical implementation (this is important for high-performance computing and crucial when we consider using quantum computers).
	
	$^{***}$In the terms of compressive sensing, the thresholded sample $T_tSf$ plays the role of the signal, usually assumed sparse, the sketching map $\Phi$ that of the measurement map, and the sketch vector $y$ of generalized moments that of the measurements. Unlike in compressive sensing, $f$ does not need to be sparse, and the decoding does not aim at recovering $T_tSf$, which would overfit the sample, but only at locating the maximum of $f$.

	In the Numerical Results, we use matching pursuit for the decoding procedure and several sketch functions to compute generalized moments. The sketch functions used in this paper are:\\
	
	\noindent {\bf Quadruplets interactions}\\
	Each line of the sketching matrix is defined by
	\begin{align*}
		\Phi^{i,i+1,i+2,i+3}_{x_1,x_2,x_3,x_4} &=1_1\otimes … \otimes 1_{i-1}\\
		&\otimes (1-x_1,x_1)\otimes(1-x_2,x_2)\\
		&\otimes (1-x_3,x_3)\otimes(1-x_4,x_4)\\
		&\otimes 1_{i+5}\otimes…\otimes 1_N
	\end{align*}
	for every starting position $1 \leq i \leq N-3$ and $x_i$s take all binary values.\\
	
	\noindent {\bf Quintuplets interactions}\\
	Each line of the sketching matrix is defined by
	\begin{align*}
		\Phi^{i,i+1,i+2,i+3,i+4}_{x_1,x_2,x_3,x_4,x_5} &=1_1\otimes … \otimes 1_{i-1}\\
		&\otimes (1-x_1,x_1)\otimes(1-x_2,x_2)\\
		&\otimes (1-x_3,x_3)\otimes(1-x_4,x_4)\\
		&\otimes(1-x_5,x_5) \otimes 1_{i+6}\\
		&\otimes…\otimes 1_N
	\end{align*}
	for every starting position $1 \leq i \leq N-4$ and $x_i$s take all binary values.\\
	
	\noindent{\bf Random sketch } A random matrix whose coefficients are drawn from a Gaussian distribution.\\
	\noindent Later in the Theoretical Analysis, we use the matching pursuit algorithm for the decoding procedure. This defines the following problem for the first crucial step of matching pursuit:\\
	
	\noindent{\bf Problem II}$\qquad\qquad x^*_{II}=\argmax \Phi^T\Phi TSf$.\\
	
	\subsection*{The Class of Compressible Combinatorial Problems}
	\label{sec:problem_setting}
	
	The MCCO algorithm described above is not assumed to solve arbitrary combinatorial problems. However, the scope of its applicability is large and contains classes of natural problems. In this section, we identify a class of objective functions that we call {\it compressible}. In the Numerical Results, we show that MCCO performs well on this class, in particular when only a few queries to the objective are available, and compare it with digital annealing.
	
	Consider a compressible objective function $f_{\mathcal{R}} : x \mapsto f_{\mathcal{R}}(x)$, where $x\in\{0,1\}^N$ and $f_{\mathcal R}$ is uniquely determined by a set $\mathcal R$ of a few rules. A rule $(r,\omega_r)$ is a pattern $r \in \{0,1,*\}^{k_r}$ of length $k_r \leq N$ with an associated reward $\omega_r \in \mathbb{R}$, where the wildcard $*$ matches either bit value. Each time the pattern $r$ matches a sub-string of $k_r$ consecutive bits of the input $x$, the total cost $f_{\mathcal{R}}(x)$ increases by $\omega_r$. A rule without wildcards constrains $k_r$ consecutive bits; we call it \emph{local}. A rule with wildcards couples bits up to $k_r - 1$ positions apart without constraining the bits in between; we call it \emph{non-local}.
	
	Therefore, for a given input $x$, the function is defined as:
	\begin{equation}
		f_\mathcal{R}(x)=\sum_{r \in \mathcal{R}}\sum_{i=1}^{N-k_r+1} \omega_r\, \delta_r(x_{i \rightarrow i+k_r-1})
	\end{equation}
	where $x_{i \rightarrow i+k_r-1}$ is the sub-string of $x$ made of the bits $i$ to $i+k_r-1$, and $\delta_r(x_{i \rightarrow i+k_r-1}) = 1$ if $x_{i+j-1} = r_j$ for every position $j$ such that $r_j \neq *$, and $0$ otherwise.
	
	The function is compressible in the sense that all the information is contained in a few values and there exists a (non-linear) mapping $F_N : \bigcup_{k=1}^{N}\{0,1,*\}^{k} \rightarrow \mathbb{R}^{2^N}$, defined on the set of patterns, embedding the rules into the optimization space. This mapping is given by the following transform:
	\begin{equation}
		\label{eq:transform}
		F_N(r)=\sum_{x\in\lbrace 0, 1\rbrace^N}\sum_{i=1}^{N-k_r+1} \delta_r(x_{i \rightarrow i+k_r-1}) \bf{x}
	\end{equation}
	where $\bf{x}$ is a vector of which the $x$-th element is $1$ and the rest is $0$. The transform $F_N$ maps a rule $r$ to a real vector of size $2^N$ that contains the values taken by the function $f_\mathcal{R}$ on $\{0,1\}^N$. We can define the following vector:
	\begin{equation}
		f_\mathcal{R} = \sum_{r \in \mathcal{R}} \omega_r F_N(r),
	\end{equation}
	as a linear combination of the transforms of all rules in $\mathcal{R}$ and thus redefine the function $f_\mathcal{R}(x)$ as the $x$-entry of the vector $f_\mathcal{R}$.\\
	
	The function $f_\mathcal{R}$ can also be written as a diagonal Hamiltonian, $\mathcal{H} = \text{diag}(f_\mathcal{R}) = \sum_{r \in \mathcal{R}} \omega_r\, \text{diag}(F_N(r))$. Using the conventional mapping from binary models to Ising Hamiltonians, we have
	\begin{equation}
		\mathcal{H} = \sum_{r \in \mathcal{R}} \sum_{i=1}^{N-k_r+1} \omega_r \bigotimes_{j \,:\, r_j \neq *} \frac{\mathbb{I}+(-1)^{r_j}\sigma^z_{i+j-1}}{2},
	\end{equation}
	where the tensor product implicitly contains the identities on the spins not constrained by the rule. Expanding the products gives a generalized Ising Hamiltonian with interactions among up to $k_r$ spins within a window of $k_r$ consecutive spins. For more details on the connections with Ising Hamiltonians, refer to the Discussion.
	
	This structure does not make the optimization of $f_\mathcal{R}$ hard in itself: when the rules are known, every rule acts on at most $k = \max_r k_r$ consecutive bits, and dynamic programming along the bit-string, with the last $k-1$ bits as state, finds the maximum exactly at a cost linear in $N$ and exponential only in $k$. The setting of this paper is different. The solver only has black-box access to $f_\mathcal{R}$: it can query it, but it is given neither the rules, nor their number, lengths or patterns. Methods that exploit the structure of the objective then cannot be applied directly, and the problem is typically approached with black-box heuristics such as simulated annealing or digital annealing. MCCO belongs to this category: it does not use the rules, but it exploits the compressibility of the objective. Objectives of this kind arise when the function is given by a simulator, an experiment, a learned model or the returns of a reinforcement learning environment (see the Discussion): their values may be built from structured rules, but the rules themselves are not available to the solver.
	
	\section*{Numerical Results}
	\label{sec:numerics}

	In this section, we test MCCO on randomly drawn compressible functions and compare it with digital annealing \cite{aramon_physics-inspired_2019}, a black-box heuristic that, like MCCO, accesses $f_\mathcal{R}$ only through evaluations.
	The methods are compared at an equal number of queries to $f_\mathcal{R}$, which is the relevant cost when $f_\mathcal{R}$ is expensive to evaluate. Our focus is the regime of few queries, where this cost matters most. All simulations use the TrOMA library \cite{baptiste_baptistechevtroma_2026}; the code, the seeds and the records of every run are available and fully reproducible (see Data Availability).\\

\begin{figure}[ht]
		\centering
		\includegraphics[width=0.8\linewidth]{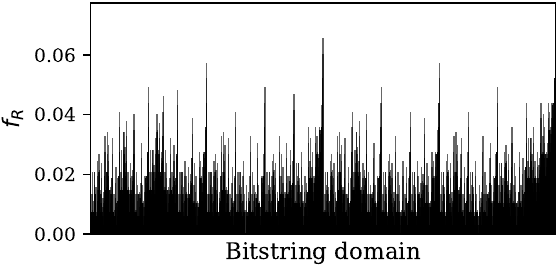}
		\captionsetup{justification=raggedright}
		\caption{Example of $f_\mathcal{R}$. On the $x$ axis, we list all binary sequences. Their respective rewards are shown on the $y$ axis.}
		\label{fig:spectrum}
	\end{figure}

	\textbf{Setup}: We take bit-strings of length $N = 20$, so that the search space contains $2^{20}$ strings and the optimum $f_\mathcal{R}(x^*)$ of every instance is known by exhaustive enumeration. The instances are drawn from the class of compressible functions $f_\mathcal{R}$ described in the Methods (an example is shown in Figure \ref{fig:spectrum}). Each rule $r \in \mathcal{R}$ has a length $k_r$ drawn uniformly in $\{4, 5, 6\}$ and a reward $\omega_r$ drawn uniformly in $(0, 1]$. We consider two families of rules: \emph{local} and \emph{non-local}. For the \emph{local} family, the patterns are drawn uniformly in $\{0,1\}^{k_r}$. For the \emph{non-local} family, a pattern has its first and last symbols drawn in $\{0,1\}$ while the interior symbols are chosen within $\{0,1,*\}$, where $*$ matches either bit value. Such a rule couples bits up to $k_r - 1$ positions apart without constraining the bits in between. For each family and each number of rules $|\mathcal{R}| \in \{1, \dots, 5\}$, we draw $100$ instances, i.e. $1000$ instances in total. The solvers are told neither the family, nor the rules, nor their lengths. Degenerate maxima are allowed, and a run is counted as a success when it returns any maximizer, i.e. $f_\mathcal{R}(\hat{x}) = f_\mathcal{R}(x^*)$.

	A run of MCCO draws $n$ strings uniformly with replacement. The budgets $n$ range from $50$ to $102400$ queries (about $10\%$ of the search space), spaced logarithmically. Every instance is attempted $10$ times at every budget, and each run draws its own sample of $102400$ strings. Then the run at budget $n$ uses its first $n$ strings, while the same sample is shared by all sketch functions. MCCO uses the same settings on all instances, without any tuning per family or per number of rules. The threshold $t$ is chosen as the $85$th percentile of the $n$ sampled values, so that it is set from the sample itself. We compare the three sketch functions of the Methods: Quadruplets, Quintuplets, and a Gaussian random sketch, drawn once per instance. The decoding procedure is matching pursuit with $5$ iterations, which returns $5$ candidate strings. We make a slight, low-cost improvement to the last step of Algorithm 1: instead of taking the largest value of $\tilde{f}$, we evaluate the $5$ candidates through $f_\mathcal{R}$, at a cost of at most $5$ additional queries, and the estimate $\hat{x}$ is the best string among the candidates and the sampled strings, whose values are already known. With the Quadruplets and Quintuplets sketches, each decoding step is the maximization of a one-dimensional Ising-type chain of $k$-local terms, solved exactly by dynamic programming along the chain, at a cost that scales exponentially in the interaction size $k$ rather than in $N$. The random sketch has no such structure and requires a search over all $2^N$ strings to find an exact solution. As it is slower to decode, we run it only on $10$ instances per family and $|\mathcal{R}|$.

	The baseline is digital annealing (DA) \cite{aramon_physics-inspired_2019}, a variant of simulated annealing for binary optimization. At each step, DA evaluates the $N$ strings at Hamming distance one from the current string and moves to one of the neighbors accepted by the Metropolis criterion, so that a budget of $n$ queries allows about $n/N$ steps; it returns the best string it has queried. Like MCCO, it uses the same setting on every instance: the temperature decreases geometrically from $T_0 = 0.36$ to $T_\mathrm{end} = 0.11$, and an offset, increased by $0.1$ after each step without an accepted move, lets it escape local maxima.

	For every run, we record the success and the distance to the optimum $(f_\mathcal{R}(x^*) - f_\mathcal{R}(\hat{x}))/\sigma_{f_\mathcal{R}}$, where $\sigma_{f_\mathcal{R}}$ is the standard deviation of $f_\mathcal{R}$ over $\{0,1\}^N$. Success rates are given with $95\%$ bootstrap intervals. Each of $2000$ bootstrap replicates draws the instances of an ensemble with replacement, then, for each drawn instance, its $10$ runs with replacement, and computes the success rate; the interval spans the $2.5$th to $97.5$th percentiles of the $2000$ replicates. Every random draw (instance, sample, random sketch, annealing run) has its own seed, derived from a single master seed. The simulations ran on an Intel Xeon Silver 4410Y processor with $48$ single-threaded worker processes.\\

\begin{figure*}[ht]
		\centering
		\captionsetup{justification=raggedright}
		\begin{subfigure}{\linewidth}
			\centering
			\includegraphics[width=0.85\linewidth]{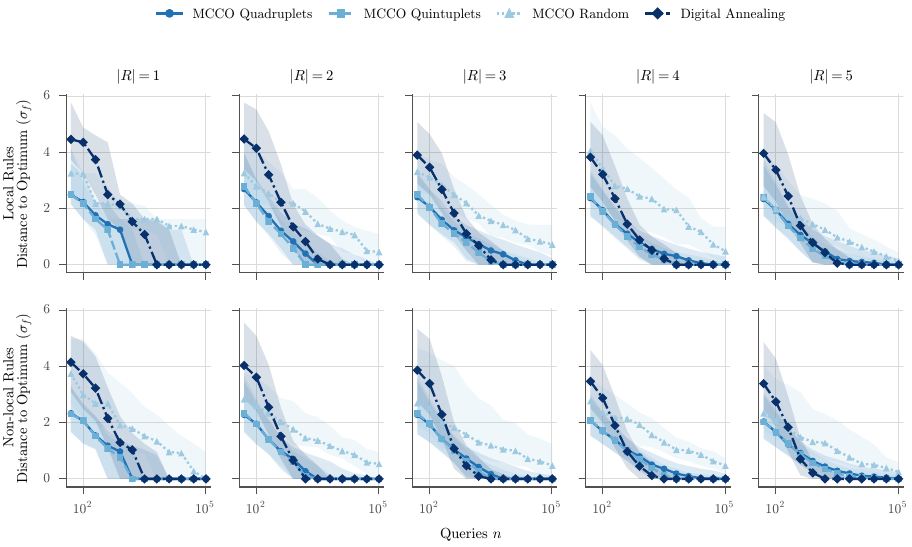}
			\caption{Distance to the optimum $(f_\mathcal{R}(x^*) - f_\mathcal{R}(\hat{x}))/\sigma_{f_\mathcal{R}}$, where $\sigma_{f_\mathcal{R}}$ is the standard deviation of $f_\mathcal{R}$ over $\{0,1\}^N$. Lines show the median over all runs, shaded bands the interquartile range.}
			\label{fig:e1_distance}
		\end{subfigure}

		\vspace{0.15cm}
		\begin{subfigure}{\linewidth}
			\centering
			\includegraphics[width=0.85\linewidth, trim={0 0 0 22pt}, clip]{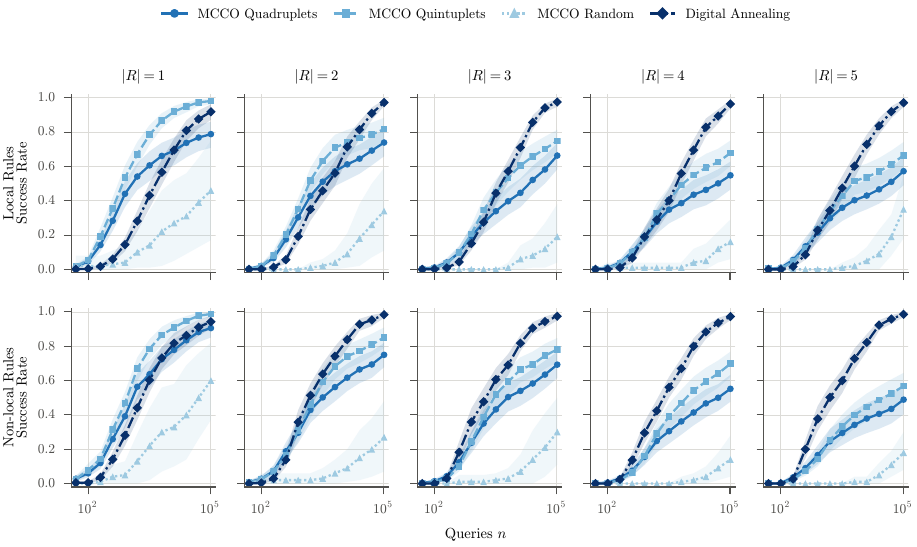}
			\caption{Success rate, $f_\mathcal{R}(\hat{x}) = f_\mathcal{R}(x^*)$. Shaded bands are $95\%$ bootstrap intervals (instances, then runs within instances).}
			\label{fig:e1_success}
		\end{subfigure}
		\caption{MCCO and digital annealing as a function of the number of queries $n$, for local rules (top row of each panel) and non-local rules (bottom row), with $|\mathcal{R}| = 1, \dots, 5$ rules and $N = 20$; $100$ instances per panel and $10$ runs per instance ($10$ instances for MCCO Random). The MCCO estimate $\hat{x}$ is the best of the sampled strings and of the $5$ matching-pursuit candidates. The largest budget, $n = 102400$, is about $10\%$ of $2^{20}$.}
		\label{fig:e1}
	\end{figure*}

	\textbf{Results interpretation}: Figure \ref{fig:e1_distance} shows the distance to the optimum as a function of the number of queries. For each method and ensemble, the line is the median over all runs and the band the interquartile range. The main result of this figure is that MCCO is efficient with few queries. Up to a few hundred queries, MCCO with the structured sketches returns better solutions than digital annealing, for both families and every number of rules. At $n = 50$, the median distance is $2.0$--$2.8\,\sigma_{f_\mathcal{R}}$ for Quadruplets and Quintuplets against $3.4$--$4.5\,\sigma_{f_\mathcal{R}}$ for digital annealing, and at $n = 200$ it is $1.2$--$1.8\,\sigma_{f_\mathcal{R}}$ against $1.8$--$3.7\,\sigma_{f_\mathcal{R}}$. These budgets correspond to the regime in which evaluations of $f_\mathcal{R}$ are the limiting resource. The difference comes from how the queries are used: digital annealing spends $N = 20$ queries per step on the neighbors of a single string, so that it can only make a few steps, while MCCO uses every sample to estimate the moments of $f_\mathcal{R}$ over the whole search space. The success rates (Figure \ref{fig:e1_success}) follow the same trend: at $n = 200$ and $|\mathcal{R}| \leq 2$, MCCO with Quintuplets finds the optimum in $7$--$19\%$ of the runs, against $1$--$4\%$ for digital annealing.

	The advantage of MCCO is specific to the low-query regime: digital annealing catches up as the budget grows. Its median distance reaches that of the best structured sketch at $n \approx 400$--$1600$ on the non-local family with $|\mathcal{R}| \geq 2$, and at $n \approx 3200$--$6400$ on the local family. In specific cases, in particular when the number of rules is small, MCCO remains ahead at every budget. For $|\mathcal{R}| = 1$, MCCO with Quintuplets has a success rate of $0.98$ (local) and $0.99$ (non-local) at the highest budget ($n = 102400$), against $0.92$ and $0.94$ for digital annealing. 

	On the other hand, for $|\mathcal{R}| \geq 2$, the success of MCCO levels off, at $0.57$--$0.85$ for Quintuplets and $0.49$--$0.75$ for Quadruplets at the highest budget. This plateau does not contradict the Theoretical Analysis, whose guarantee holds only when the sketch preserves the maximum. The success rate also appears correlated with the fraction of instances on which the sketch preserves the maximum. However, this condition is only sufficient, and some instances are solved without it, so the correlation does not explain the plateau. 

	Figure \ref{fig:e1} also shows that the structure of the sketch matrix is a key element of MCCO, with the Quadruplets and Quintuplets sketches achieving good performance in every scenario, while the random sketch struggles to propose any useful candidate.

	The results on the local and non-local families are similar, with the same settings and without being told which family an instance belongs to. MCCO does not require the sketch to match the rules either: on an instance made of $3$ rules of length $4$, the matched sketch (Quadruplets) and a mismatched one (Quintuplets) reach comparable success rates, e.g. $0.31$ and $0.33$ at $n = 3200$ and $0.41$ and $0.37$ at $n = 102400$ over $300$ runs (see Supplementary Figure 1).

	Figure \ref{fig:e2_threshold} isolates the effect of the threshold on one instance of each family with $|\mathcal{R}| = 5$, at $n = 6400$, with $300$ runs per threshold and the matching-pursuit candidates as estimate. Thresholding increases the success probability by a factor of $2$ to $5$ with respect to $t = 0$: from $0.11$ to $0.45$ (Quadruplets) and from $0.34$ to $0.71$ (Quintuplets) on the local instance, and from $0.13$ to $0.47$ and from $0.12$ to $0.61$ on the non-local instance. The success probability varies little for thresholds between the $50$th and the $95$th percentiles of $f_\mathcal{R}$, a range that contains the $85$th percentile used in the other experiments, so that MCCO is not sensitive to the precise choice of $t$. Above the $95$th percentile, the success drops ($0.02$--$0.10$ at the $99.9$th percentile), as only a few samples remain above the threshold and the moments of $T_t f_\mathcal{R}$ are poorly estimated. This agrees with the guarantee of Theorem \ref{theo:threshold} below.

	Table \ref{tab:e4_cost} gives the wall-clock time of each stage at a budget of $n = 6400$. With the structured sketches, a complete MCCO run takes $0.18$--$0.19$\,s, of which $0.15$\,s is spent in the sampling stage, while sketching and decoding take $0.01$--$0.02$\,s each. Decoding the random sketch takes $11$--$11.5$\,s, about $1000$ times longer, because of the search over the $2^N$ strings mentioned above. Digital annealing takes $0.007$\,s. In these simulations, a query to $f_\mathcal{R}$ is a table look-up, so the times reflect the implementations rather than the cost of the queries. When $f_\mathcal{R}$ is costly to evaluate, as in the black-box setting, the cost of both methods is set by the number of queries, which is the axis of Figure \ref{fig:e1}, and MCCO is ahead at small budgets.

	All in all, MCCO finds good solutions to compressible problems from a small number of queries, in the regime where the number of evaluations of $f_\mathcal{R}$ is the limiting resource, and it does so with the same settings on both families of rules, without knowing their structure. At larger budgets, digital annealing reaches a higher success rate when the problem contains several rules. The open-source TrOMA library \cite{baptiste_baptistechevtroma_2026} allows anyone to use MCCO and to reproduce our results. As the scope of problems our method handles best is not yet perfectly known, we invite the reader to try our method on their own problems through the library.

	\section*{Theoretical Analysis}
	\label{sec:guarantees_new}
	
	In this section, we give theoretical evidence of the soundness of the MCCO method. In other words, we show that the solution of Problem II is consistent with, or approximately matches, that of Problem I. 
	
	\subsection*{Matrices that preserve the maximum}
	
	Let us discuss the properties of appropriate matrices $\Phi$ to be used as sketch maps in MCCO (Algorithm 1). 
	It is possible to characterize good sketch matrices depending on the parameters of the problem, e.g. the gap between the two largest values of the tested function. Although a definition based on quantitative assumptions remains an open problem, our numerical results demonstrate the existence of suitable sketch matrices. In addition, we give qualitative properties for matrices belonging to this class. 
	
		\begin{definition}[Surrogate and empirical surrogate]
		\label{def:surrogate}
		Let $G = \Phi^T \Phi$ be the Gram matrix. Let $g^x$ denote the $x$-th column of $G$, and, for $s \in \lbrace 0,1 \rbrace^N$, let $g^x(s)$ and $f(s)$ denote the values taken at $s$ by that column and by $f$ respectively. The \emph{surrogate function} associated with $f$ and $\Phi$ is
		\begin{equation*}
			F(x) = 2^{-N}(Gf)_x = \mathbb{E}_s\left[f(s)g^x(s)\right], \qquad x \in \lbrace 0,1 \rbrace^N,
		\end{equation*}
		where the expectation is taken over $s$ uniformly distributed on $\lbrace 0,1 \rbrace^N$. Given a sample $S_n = (s_1,\ldots,s_n)$ drawn i.i.d. uniformly on $\lbrace 0,1 \rbrace^N$ with replacement, the \emph{empirical surrogate} is
		\begin{equation*}
			\tilde{F}_{S_n}(x) = \frac{1}{n}\sum_{k=1}^{n} f(s_k)g^x(s_k).
		\end{equation*}
	\end{definition}
	
	\noindent By construction $\mathbb{E}\left[\tilde{F}_{S_n}(x)\right] = F(x)$ for every $x$. Since $F$ differs from $Gf$ only by the positive factor $2^{-N}$, the two have the same maximizer, and Problem II is the maximization of an empirical surrogate of this kind, built from the sampled and thresholded values of $f$.
	
	We consider under which conditions the maximum of the objective function $f$ coincides with the maximum of the surrogate $F(x)$. Then, we will give evidence that the maximum can survive under sampling, i.e. $F(x)$ and $\tilde{F}_{S_n}(x)$ have the same maximum.
	First, let us consider some general properties of $G$:   
	
	\begin{itemize}
		\item $G$ is a square matrix whose dimensions are $2^N \times 2^N$
		\item $G$ is symmetric, i.e. $G = G^T$
		\item $\Phi$ can be normalized such that $G$ only has ones on its diagonal
	\end{itemize}

	As an additional desirable property, $G$ should behave as a quasi-isometry. For example, this is the case for matrices satisfying the RIP, a well-studied property in compressive sensing, or more generally for matrices satisfying the Johnson-Lindenstrauss lemma. Notice that, just like the identity matrix, $G$ already has ones on its diagonal. However, we can allow additional degrees of freedom without changing the position of the maximum of $f$. In order to achieve quasi-isometry, the maximum of each row should remain on $G$'s diagonal, and all off-diagonal elements should be relatively small. These are already the assumed properties for the measurement matrices used in compressive sensing. The off-diagonal elements of $G$ are described by the coherence $\mu$ of $\Phi$ and, in fact, since diagonal elements of $G$ are equal to one, we want $\mu$ smaller than $1$. However, our numerical results show that too low coherence gives a weak generalization power, which is an issue when we later introduce sampling.
	
	For a structured matrix with nearest-neighbor $k$-uplet patterns like those defined after the algorithm description, the coherence is known and grows as $\mu = \frac{N-k}{N-k+1}$ (see Supplementary Information, Section 2 for the proof). The quintuplets structured matrix used in the Numerical Results has coherence $\mu = \frac{15}{16}$, and the random matrices used in the same section have a coherence of about $0.31$. 
	
	Other quantitative assumptions on $G$ include that any two rows of $G$ should not be too similar or too different. In fact, problems $f$ whose structure contains a gap between the two largest values should typically be solved under this assumption.
	Notice that because $\Phi$ contains only $m$ rows, the rank of $G$ is smaller than or equal to $m$. Thus, $G$ is not full-rank. Because of the presence of eigenspaces that are mapped to zero, without the above assumption, $G$ tends to close the gap in $f$. If the rows are too similar, inputs are easily mixed; if they are too different, we amplify the chance to promote a wrong candidate for the maximum. The gap assumption mentioned above is common in optimization algorithms. It is a requirement of many meta-heuristic methods, for instance quantum annealing, where one has to consider the gap in the optimized Hamiltonian.
	
	An additional assumption on $G$ is that it should be easily decodable. This assumption is not needed for $G$ to preserve the maximum. However, in order for the surrogate $F$ to be a simpler optimization problem than $f$, this property is highly desirable. With that goal in mind, we can artificially limit the choice of $G$ to a class that implies a tractable problem. The structured matrices that we suggested in the Methods implement the feasibility assumption. Indeed, the Quadruplet and Quintuplet matrices represent a spin-chain Ising problem and can be solved by dynamic programming, e.g., \cite{schuch_matrix_2010}. In our previous work \cite{chevalier_compressed_2024}, we also discussed how more complex structured matrices can be decoded using quantum computers (see the Discussion for more details).
	
	\subsection*{Preserving the maximum under sampling}
	
	In the previous section, we discussed the existence of sketch matrices which do not change the maximum of the sketched function. In this section, we study how these matrices behave when values of the objective function $f$ are being sampled. Under large enough sampling, we quantify the probability that the maximizer of the empirical surrogate $\tilde{F}_{S_n}$ of Definition \ref{def:surrogate} coincides with the maximizer of the surrogate $F$. We also show that we can benefit from applying the thresholding operator $T_t$ in this task.
	
	The following theorem bounds the probability that the maximizer of the surrogate is displaced when the surrogate is estimated from a finite sample.
	
	\begin{theorem}
		\label{theo:max_preserve_samp}
		Let $x^*$ be the position of the maximum of $f$, and assume that $G$ preserves this maximum when applied to the full function, i.e. $F(x^*) > F(x)$ for all $x \neq x^*$. For each $x \neq x^*$, we define the gap random variable $\Delta_x$, its average $\Theta_x$ and its second moment $\nu_x^2$
		\begin{equation*}
			\Delta_x(s) = f(s)\left(g^{x^*}(s)-g^x(s)\right), \qquad \Theta_x = \mathbb{E}_s\left[\Delta_x(s)\right] = F(x^*)-F(x) > 0, \qquad \nu_x^2 = \mathbb{E}_s\left[\Delta_x(s)^2\right].
		\end{equation*}
		We also define
		\begin{equation*}
			\Theta_{\min} = \min_{x \neq x^*} \Theta_x = F(x^*) - \max_{x\neq x^*} F(x), \qquad \nu^2 = \max_{x \neq x^*} \nu_x^2, \qquad M = \max_{x \neq x^*} \left\| \Delta_x - \Theta_x \right\|_\infty.
		\end{equation*}
		The following property holds:
		\begin{equation}
		\label{eq:sampling_max}
			P\left(\argmax_x \tilde{F}_{S_n}(x) \neq x^* \right) \leq 2^N \exp\left(\frac{-n\Theta_{\min}^2}{2\nu^2+\frac{2}{3}M\Theta_{\min}}\right),
		\end{equation}
		and in particular the maximizer is preserved with probability at least $1-\delta$ whenever
		\begin{equation}
		\label{eq:sample_size}
			n \geq \left(\frac{2\nu^2}{\Theta_{\min}^2} + \frac{2M}{3\Theta_{\min}}\right)\left(N\ln 2 + \ln(1/\delta)\right).
		\end{equation}
		Both bounds remain valid when $M$ is replaced by any upper bound on it.
	\end{theorem}
	
	\begin{proof}
		The domain $\lbrace 0,1 \rbrace^N$ is finite and $f$ is bounded, so each $\Delta_x$ is a bounded random variable under uniform sampling. Since the $s_k$ are drawn independently, the variables $\Delta_x(s_1),\ldots,\Delta_x(s_n)$ are i.i.d. with common mean $\Theta_x \geq \Theta_{\min} > 0$ and variance at most $\nu_x^2 \leq \nu^2$. The constant $M$ is finite for the same reason, and admits the explicit bound
		\begin{equation*}
			M \leq 2\|f\|_\infty \max_{x\neq x^*,\, s}\left|g^{x^*}(s)-g^x(s)\right| \leq 2m\|f\|_\infty,
		\end{equation*}
		where $m$ is the number of rows of $\Phi$, since $G=\Phi^T\Phi$ has non-negative integer entries bounded by $m$.\\
		For every $x \neq x^*$, we define the summed random variable $\tilde{\Delta}_x$ associated with a sample set $S_n$:
		\begin{equation*}
			\tilde{F}_{S_n}(x^*)-\tilde{F}_{S_n}(x) = \frac{1}{n}\sum_{k=1}^{n}\Delta_x(s_k) =: \tilde{\Delta}_x .
		\end{equation*}
		Hence $x^*$ is the unique maximizer of $\tilde{F}_{S_n}$ if and only if $\tilde{\Delta}_x > 0$ for all $x \neq x^*$. The equivalence holds since the gap exists for all competitors.\\
		Consequently, using the union bound,
		\begin{equation*}
			P\left(\argmax_x \tilde{F}_{S_n}(x) \neq x^*\right) = P\left(\bigcup_{x \neq x^*}\lbrace \tilde{\Delta}_x \leq 0\rbrace\right) \leq \sum_{x\neq x^*} P\left(\tilde{\Delta}_x \leq 0\right).
		\end{equation*}
		For i.i.d. variables with $\|\Delta_x - \Theta_x\|_\infty \leq M$ and variance at most $\nu_x^2$, Bernstein's concentration inequality (whose bound increases with the variance and with $M$) gives, for all $t>0$,
		\begin{equation*}
			P\left(\tilde{\Delta}_x - \Theta_x \leq -t\right) \leq \exp\left(\frac{-nt^2}{2\nu_x^2+\frac{2}{3}Mt}\right).
		\end{equation*}
		Choosing $t = \Theta_x > 0$, the left-hand side becomes $P(\tilde{\Delta}_x \leq 0)$, so that
		\begin{equation*}
			P\left(\tilde{\Delta}_x \leq 0\right) \leq \exp\left(\frac{-n\Theta_x^2}{2\nu_x^2+\frac{2}{3}M\Theta_x}\right).
		\end{equation*}
		Finally, the map $u \mapsto u^2/(2\nu^2+\frac{2}{3}Mu)$ is increasing on $u>0$, so $\Theta_x \geq \Theta_{\min}$ and $\nu_x^2 \leq \nu^2$ yield
		\begin{equation*}
			P\left(\tilde{\Delta}_x \leq 0\right) \leq \exp\left(\frac{-n\Theta_{\min}^2}{2\nu^2+\frac{2}{3}M\Theta_{\min}}\right).
		\end{equation*}
		Summing over the $2^N-1 < 2^N$ competitors results in equation (\ref{eq:sampling_max}). Requiring the resulting bound to be at most $\delta$ and solving for $n$ gives equation (\ref{eq:sample_size}).
	\end{proof}
	
	Theorem \ref{theo:max_preserve_samp} characterizes the behavior of the algorithm under sampling. The probability of success is controlled by the mean, the second moment and the range of the random variables $\Delta_x$, and in particular by $\Theta_{\min}$, the gap of the surrogate function. This becomes a statement about the objective itself under the assumptions made above on $G$: when the sketch acts as a quasi-isometry in the restricted sense of compressive sensing, the surrogate inherits the shape of $f$, and $\Theta_{\min}$ reflects the gap between the two largest values of the objective function. The dimension of the problem then enters the required sample size only through a factor $N\ln 2$, i.e. logarithmically in the size of the search space, the remaining dependence being carried by $\Theta_{\min}$, $\nu^2$ and $M$. We tested Theorem \ref{theo:max_preserve_samp} numerically on instances of the Numerical Results; the results are given in the Supplementary Information, Section 3.
	
	\subsection*{Effect of Thresholding}
	
	In this subsection, we show that the guarantee of Theorem \ref{theo:max_preserve_samp} still holds when the sample is thresholded, and then that thresholding improves it.	
	The thresholding step changes the sampled function from $f$ to $T_t f$, where $$T_t f(x) = \begin{cases}	f(x) & \text{if } f(x) \geq t \\  0 & \text{otherwise} \end{cases}.$$ 
	Let $f : \{0,1\}^N \rightarrow \mathbb{R}^+$ have a unique maximizer $x^*$ and let $x_2 = \argmax_{x \neq x^*} f(x)$ be its second maximum.
	The following proposition shows that, for an admissible threshold, this change leaves both the maximizer and the gap of $f$ untouched.
	
	\begin{proposition}
		\label{prop:threshold}
		Let $\Delta = f(x^*) - f(x_2) > 0$ be the gap of $f$. If $0 \leq t \leq f(x_2)$, then $x^*$ is the unique maximizer of $T_t f$ and the gap of $T_t f$ is exactly $\Delta$.
	\end{proposition}
	\begin{proof}
		The proof is given in the Supplementary Information, Section 4.
	\end{proof}

	Since $T_t$ acts pointwise, for a fixed $t$ it commutes with the sampling operator: $T_t S f = S T_t f$. Thresholding the sample is therefore the same as sampling the thresholded function, so the results on sampling apply to $T_t f$ in place of $f$.

	In this analysis, the threshold $t$ is considered fixed, whereas in practice it is chosen from the sampled values. 
	One can consider that a first sample is used to decide $t$, within the admissible range $0 \leq t \leq f(x_2)$. 
	For instance, in the Numerical Results section, the threshold $t$ was set to the 85th percentile of the sampled values.
	Once the value of $t$ is fixed, a new sample, independent of the first one, is drawn, so that the resulting empirical surrogate is unbiased for the thresholded function $T_t f$. 
	This is only to make the analysis rigorous; in practice, the samples used before and after fixing $t$ are similar, so the same sample can be used throughout.
	
	\begin{definition}[Thresholded quantities]
		\label{def:thresholded}
		Let $F_t$ be the surrogate (Definition \ref{def:surrogate}) of $T_t f$. For $x \neq x^*$, the quantities of Theorem \ref{theo:max_preserve_samp} computed from $T_t f$ instead of $f$ are
		\begin{equation*}
			\Delta_{x,t}(s) = T_t f(s)\left(g^{x^*}(s)-g^x(s)\right), \qquad \Theta_{x,t} = \mathbb{E}_s\left[\Delta_{x,t}(s)\right] = F_t(x^*)-F_t(x), \qquad \nu_{x,t}^2 = \mathbb{E}_s\left[\Delta_{x,t}(s)^2\right],
		\end{equation*}
		and we let $\Theta_{\min,t} = \min_{x \neq x^*}\Theta_{x,t}$, $\nu_t^2 = \max_{x \neq x^*}\nu_{x,t}^2$ and $M_t = \max_{x \neq x^*}\left\|\Delta_{x,t} - \Theta_{x,t}\right\|_\infty$. We also let $(T_t f)^\complement = f - T_t f$ be the part of $f$ discarded by thresholding, and $F^\complement_t$ its surrogate.
	\end{definition}

	\begin{corollary}[of Theorem \ref{theo:max_preserve_samp}, Thresholded Function]
		\label{cor:threshold}
		Let $0 \leq t \leq f(x_2)$, assume that $G$ preserves the maximum of $T_t f$, and let $\tilde{F}_{S_n,t}(x) = \frac{1}{n}\sum_{k=1}^{n} T_t f(s_k)g^x(s_k)$ be the empirical surrogate of $T_t f$. Then, with $x^*$ the maximizer of $f$,
		\begin{equation}
		\label{eq:sampling_max_t}
			P\left(\argmax_x \tilde{F}_{S_n,t}(x) \neq x^* \right) \leq 2^N \exp\left(\frac{-n\Theta_{\min,t}^2}{2\nu_t^2+\frac{2}{3}M_t\Theta_{\min,t}}\right).
		\end{equation}
	\end{corollary}
	\begin{proof}
		Theorem \ref{theo:max_preserve_samp} only requires the sampled function to be bounded, so it applies to $T_t f$, whose maximizer is $x^*$ by 
		Proposition \ref{prop:threshold}.
	\end{proof}

	Thresholding acts on both the second moment $\nu_t^2$ and the surrogate gap $\Theta_{\min,t}$ in the bound (\ref{eq:sampling_max_t}). 
	The second moment $\nu_t^2$ decreases as $t$ increases, and the surrogate gap $\Theta_{\min,t}$ is at least $\Theta_{\min}$ whenever $x^*$ is a 
	minimizer of $F^\complement_t$, the surrogate of the discarded part (Supplementary Information, Lemmas 1 and 2).

	\begin{theorem}[Thresholding improves the guarantee]
		\label{theo:threshold}
		Let $0 \leq t \leq f(x_2)$, assume that $G$ preserves the maximum of $f$ and that $x^*$ is a minimizer of $F^\complement_t$, and let $\bar{M} = 2\|f\|_\infty \max_{x\neq x^*,\, s}\left|g^{x^*}(s)-g^x(s)\right|$. Then Corollary \ref{cor:threshold} applies, and
		\begin{equation*}
			\frac{\Theta_{\min,t}^2}{2\nu_t^2+\frac{2}{3}\bar{M}\Theta_{\min,t}} \geq \frac{\Theta_{\min}^2}{2\nu^2+\frac{2}{3}\bar{M}\Theta_{\min}}.
		\end{equation*}
		Hence, with $M$ replaced by $\bar{M}$, the bound of equation (\ref{eq:sampling_max}) and the sample size of equation (\ref{eq:sample_size}) for $T_t f$ are smaller than or equal to those for $f$.
	\end{theorem}
	\begin{proof}
	The proof is given in the Supplementary Information, Section 4.
	\end{proof}

	Thresholding can be seen as rejection sampling: sampled points with $f(s) < t$ are rejected, so that the surrogate is effectively estimated from points drawn uniformly on the region $\{f \geq t\}$, which encodes the prior knowledge that the maximum lies among the large values of $f$ and that small values carry little information about its position. When $t$ becomes too large, the discarded part starts to include points that resemble $x^*$, and the gap decreases.

	In the threshold sweep of the Numerical Results, thresholding improves the success probability over a broad range of thresholds, and the success probability decreases once $t$ becomes too large.
	
	\subsection*{Specific case: Preserving the maximum with a structured sketch}
	
	We now specialize the analysis to the case where the sketch matrix is structured, e.g. chosen as Quadruplets or Quintuplets, which makes the surrogate of Definition \ref{def:surrogate} explicit. The generalized moments one estimates are then expectation values of the Walsh characters $\chi_A(x) = (-1)^{\sum_{j \in A} x_j}$, i.e. the sketch vector $y$ contains estimates $\hat{a_A}$ of the Walsh-Hadamard~\cite{odonnell_boolean_2014} coefficients $a_A = \mathbb{E}_x[f(x)\chi_A(x)]$ of $f$, indexed by the set $K$ of characters carried by $\Phi$. Written in this basis,
	\begin{equation}
		F(x) = \sum_{A \in K} a_A \chi_A(x), \qquad \tilde{F}_{S_n}(x) = \sum_{A \in K} \hat{a_A} \chi_A(x),
	\end{equation}
	so the surrogate is the truncation of the Walsh-Hadamard expansion of $f$ to $K$, and the empirical surrogate is the function actually reconstructed by the algorithm. Problem II is therefore the maximization of $\tilde{F}_{S_n}$ over $\lbrace 0,1 \rbrace^N$. The two propositions below, together with the lemma connecting them, make the resulting error budget explicit.
	
	\begin{proposition}
		\label{prop:wh_samp}
		Given a function $f$, learning $k$ of its Walsh-Hadamard coefficients, all within a given error $\epsilon$, can be achieved with sample complexity $n = O(\frac{\log{k}}{\epsilon^2})$.
	\end{proposition}
	\begin{proof}
		Each coefficient $a_A$ is an expectation of a bounded random variable. By Hoeffding's inequality, $n = O(\frac1{\epsilon^2})$ samples suffice to estimate a single coefficient within error $\epsilon$. A union bound over $k$ coefficients completes the proof.
	\end{proof}
	
	Notice that the sample complexity is independent of $N$ and remains tractable even when the size of the space becomes exponentially large.
	
	\begin{lemma}
		\label{lem:linf}
		Let $K$, with $|K| = k$, be the set of characters carried by $\Phi$. If $|a_A - \hat{a_A}| \leq \epsilon$ for every $A \in K$, then
		\begin{equation}
			\left\| f - \tilde{F}_{S_n} \right\|_\infty \leq k\epsilon + \sum_{A \notin K} |a_A| .
		\end{equation}
	\end{lemma}
	\begin{proof}
		For any $x \in \lbrace 0,1 \rbrace^N$, every character satisfies $\chi_A(x) \in \lbrace -1,1 \rbrace$, hence $\left|\chi_A(x)\right| \leq 1$. Splitting
		\begin{equation*}
			f - \tilde{F}_{S_n} = \sum_{A \in K}\left(a_A - \hat{a_A}\right)\chi_A(x) + \sum_{A \notin K} a_A \chi_A(x)
		\end{equation*}
		into the monomials of $K$, whose coefficients are estimated within $\epsilon$, and those outside $K$, which contribute their full value, and applying the triangle inequality gives the result.
	\end{proof}
	
	The two terms are the two error sources of the method: $k\epsilon$ bounds $\|F - \tilde{F}_{S_n}\|_\infty$, the sampling error already controlled by Theorem \ref{theo:max_preserve_samp}, while $\sum_{A \notin K}|a_A|$ equals $\|f - F\|_\infty$, the sketching error coming from the choice of $K$. Which coefficients are missing thus matters only through their $\ell_1$ mass, and compressibility of $f$ is precisely the assumption that this mass is small.
	
	\begin{proposition}
		\label{prop:gap}
		Given a function $f$, $f(x^*)$ its maximum and $\Delta = f(x^*) - \max_{x\neq x^*} f(x)$ the gap between the two largest values, let $\tilde{F}_{S_n}$ be the empirical surrogate of Definition \ref{def:surrogate}. If the reconstruction error $\| f - \tilde{F}_{S_n}\|_\infty < \Delta / 2$, the position of the maximum remains unchanged, i.e.
		\begin{equation}
			\argmax_x \tilde{F}_{S_n}(x) = \argmax_x f(x) = x^*.
		\end{equation}
	\end{proposition}
	\begin{proof}
		Assume $\|f - \tilde{F}_{S_n}\|_\infty < \Delta/2$. Then for any $x\neq x^*$:
		\begin{align}
			&\tilde{F}_{S_n}\left(x^*\right) \geq f\left(x^*\right)-\frac{\Delta}{2}, \label{eq:bound1}\\
			&\tilde{F}_{S_n}(x) \leq f(x)+\frac{\Delta}{2} \leq f\left(x^*\right)-\Delta+\frac{\Delta}{2}=f\left(x^*\right)-\frac{\Delta}{2}. \label{eq:bound2}
		\end{align}
		Combining (\ref{eq:bound1}) and (\ref{eq:bound2}), we have $\tilde{F}_{S_n}(x^*) \geq \tilde{F}_{S_n}(x)$ for all $x \neq x^*$. In other words, the argmax is preserved.
	\end{proof}
	
	Combining Proposition \ref{prop:wh_samp} and Proposition \ref{prop:gap} through Lemma \ref{lem:linf}, the maximum is preserved as soon as $k\epsilon + \sum_{A \notin K}|a_A| < \Delta/2$. Taking $\epsilon = \Delta/(4k)$, this holds whenever the discarded coefficients carry a mass smaller than $\Delta/4$, at a sample complexity $n = O(k^2\log k/\Delta^2)$, independent of $N$. The condition of Proposition \ref{prop:gap} is sufficient but not necessary, so the maximum can survive even larger perturbations.

	\section*{Discussion on applications and quantum computers}
	\label{sec:applications}
	
	Compressible combinatorial problems naturally show up in many fields. To see that, let us discuss learning an optimal policy in Reinforcement Learning, and finding the ground state of an Ising Hamiltonian. Reinforcement Learning is a crucial technique in today's machine learning due to its central role in several applications, such as optimal control theory, robotics, or game theory. Among the approaches to tackle it, Monte-Carlo methods solve the task based on complete returns obtained from full episodes; the name comes from the fact that one is sampling through the space of complete returns. Consider the following "game": at each step, the agent has to choose between two possible actions labeled $0$ and $1$, and a total reward is granted only after a complete episode, which consists of a sequence of $N$ actions. The goal is to find a nearly optimal policy maximizing the total reward. Since one only has access to the final returns, with the reward mechanism in each episode acting like a black box, the use of Monte-Carlo methods seems reasonable in this context. Thus, this reinforcement learning scenario exhibits a problem structure similar to the combinatorial optimization of a compressible function, and our approach can be used here as a competitive method in a context where policy networks could not be exploited to their fullest potential.
	
	Another such problem is finding the ground state of an Ising Hamiltonian. Generalized Ising models are known to be intractable for classical computers. A random-energy Ising model is given by
	\begin{equation}
		\mathcal{H}(\sigma) = \sum_{\sigma_s \in P } \omega_s \sigma_s
	\end{equation}
	where $\sigma_s$ are Pauli-Z strings within the Pauli group $P$ of given length $s$. We already made a parallel between such generalized Ising Hamiltonians and our cost function; our method is a promising new approach to tackle the problem of finding their ground state, as no efficient classical methods are known (and there will be none if $P\neq NP$) \cite{barahona_computational_1982}. 
	
	Quantum computers also seem to have a crucial role to play here, especially since so many cost functions can be mapped to Ising-like Hamiltonians. Even though they might not be able to solve the optimization problem in polynomial time (if $\text{NP} \not\subset \text{BQP}$) \cite{aaronson_np-complete_2005}, they might provide good heuristics to tackle those kinds of problems (see quantum annealing \cite{kadowaki_quantum_1998} and QAOA \cite{farhi_quantum_2014}). In our previous work \cite{chevalier_compressed_2024}, using structured maps similar to the Quadruplets and Quintuplets defined in the Methods, the optimization step in matching pursuit reduced to finding the ground state of an Ising Hamiltonian. If the sketch contains patterns on arbitrary pairs of bits, the decoding step becomes the ground-state problem of an Ising Hamiltonian on a general interaction graph, which is NP-hard \cite{barahona_computational_1982}; heuristics such as quantum annealing or QAOA can then be used for this step. More generally, when applying the Monte-Carlo Compressive Optimization algorithm with a sketch function made of structured patterns, we map the complex Ising problem into an Ising problem whose optimum is preserved. The mapping is done using a sample of the solution space, hence it is computable in tractable time.
	
	Another particularly important merit of our method is that it opens the way to optimization being carried out on specific devices such as digital annealers or quantum computers. Such devices often have limited connectivity or considerable error rates, and our method allows the model to be adapted to fit these restrictions. In two previous works \cite{seki_initialization_2024, aoki_formulation_2025}, the authors tried a similar approach by learning a surrogate model through sampling so that the optimization could be run on a quantum annealer. Extending the scope of problems that can benefit from these devices is an active area of research, and our method can be seen as a refined version of these two works, where we allow more general models and give sound theoretical evidence.
	
	\section*{Conclusion}
	\label{sec:conclusion}
	
	We introduced Monte-Carlo Compressive Optimization, a method for black-box combinatorial optimization that only queries the objective and only assumes that it is compressible. 
	MCCO samples the objective uniformly, thresholds the sample, compresses it into a sketch of generalized moments and decodes this sketch with matching pursuit, 
	without overfitting the sample itself, to propose candidates for the maximum. 
	The structure of the objective does not need to be known: it is enough that a sketch function captures the part of the objective that determines its maximum.
	
	Numerically, MCCO shows several interesting behaviors. It is most useful when queries are scarce: with a few hundred evaluations of the objective, 
	it returns better solutions than digital annealing. Every sample can contribute to an estimate of the objective over the whole search space rather 
	than just determining a local move. It is robust to agnosticity: the same settings work on local and non-local rules, and the sketch does not need to match 
	the rules of the objective. Its performance comes from the structure of the sketch rather than from the sampling alone, and it is insensitive to the 
	precise choice of the threshold. 
	
	In addition, the theoretical analysis supports these results. We gave a finite-sample bound on the probability that the empirical surrogate fails to preserve the maximum. 
	The sample size that guarantees success with high probability is only logarithmic in the size of the search space. 
	We showed that thresholding keeps the maximizer of the objective, and that, under a condition on the discarded part of the objective, it even improves the bound. 
	For structured sketches, we related the error of the method to the sampling error of the estimated coefficients and to the $\ell_1$ mass of the coefficients 
	the sketch discarded, which is where compressibility enters quantitatively.
	
	Natural next steps are to study other problem classes and larger sizes, and to choose sketch functions tailored to the problem, to the hardware, or by heuristics.
	The open-source library TrOMA makes MCCO available to test on any black-box combinatorial problem for which a cost function can be queried, and we encourage readers to take part in exploring the capacity of this method.
	
	\bibliography{ref}

	\section*{Acknowledgements}
	
	This work was supported by JST Moonshot R\&D, Grant No. JPMJMS226C and Grant No. JPMJMS256E, JST ASPIRE, Grant No. JPMJAP2427, JST COI-NEXT Grant No. JPMJPF2221, JST SPRING Grant No. JPMJSP2123.
	
	\section*{Author contributions statement}
	
	B.C. conceived the main algorithm and the numerical experiments. B.C., S.Y. and W.R. equally participated in designing the theoretical foundations and proofs. B.C. wrote the majority of the manuscript. W.R. and M.T. supervised the study and provided guidance on the ideas presented in this work. All authors reviewed the manuscript.
	
	\section*{Data Availability}
	
	The code used to produce the numerical results presented in this work is available as part of the open-source TrOMA library at https://github.com/baptistechev/TrOMA. The code used to run the numerical experiments and produce the figures is available at https://github.com/baptistechev/mcco\_numerics.
	
	\section*{Additional information}
	
	The authors declare no competing interests.
	
	\begin{figure}[ht]
		\centering
		\includegraphics[width=0.8\linewidth]{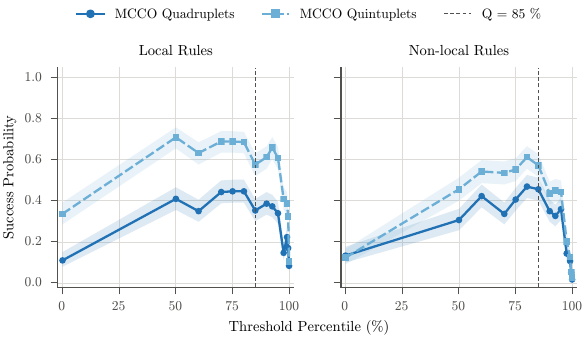}
		\captionsetup{justification=raggedright}
		\caption{Success probability as a function of the threshold $t$, given as its percentile among the $2^{20}$ values of $f_\mathcal{R}$ (the fraction of the candidates discarded, e.g. $0$: no thresholding), 
		on one instance per family with $|\mathcal{R}| = 5$ rules and a unique maximizer. Instances are chosen with a typical success probability at threshold $0$ among the instances of 
		their ensemble, at $n = 6400$. Each point is estimated from $300$ runs, with $95\%$ Wilson intervals. The estimate is the best of the $5$ matching-pursuit 
		candidates only, so that the figure shows the effect of the threshold on the decoding. The dashed line marks the $85$th percentile, the value of $q$ used 
		in Figure~\ref{fig:e1}, where $t$ is the $85$th percentile of the sampled values.}
		\label{fig:e2_threshold}
	\end{figure}
	
	\begin{table}[ht]
		\centering
		\input{fig/e4_cost.tex}
		\caption{Median wall-clock time (s) of each stage of a run at $n = 6400$, over the runs of Figure~\ref{fig:e1}, for local and non-local rules. 
		For MCCO, the stages are sampling (queries to $f_\mathcal{R}$), sketching, decoding and evaluation of the candidates; digital annealing is given as a total. 
		With the quadruplet and quintuplet sketches, $\Phi^T r$ is a one-dimensional Ising-type chain of $k$-local terms, so each decoding step $\argmax_x (\Phi^T r)(x)$ 
		is solved exactly by dynamic programming along the chain in $O((N-k+1)\,2^k)$ operations; the random sketch has no such structure and requires a search over 
		all $2^N$ strings to find the exact maximizer, hence its decoding time.}
		\label{tab:e4_cost}
	\end{table}

	\clearpage
	\setcounter{section}{0}
	\setcounter{figure}{0}
	\setcounter{table}{0}
	\setcounter{equation}{0}
	\setcounter{theorem}{0}
	\setcounter{lemma}{0}
	\setcounter{definition}{0}
	\setcounter{proposition}{0}
	\renewcommand{\figurename}{Supplementary Figure}
	\renewcommand{\theHsection}{SI.\arabic{section}}
	\renewcommand{\theHfigure}{SI.\arabic{figure}}
	\renewcommand{\theHtable}{SI.\arabic{table}}
	\renewcommand{\theHequation}{SI.\arabic{equation}}
	\renewcommand{\theHtheorem}{SI.\arabic{theorem}}
	\renewcommand{\theHlemma}{SI.\arabic{lemma}}
	\renewcommand{\theHdefinition}{SI.\arabic{definition}}
	\renewcommand{\theHproposition}{SI.\arabic{proposition}}

\section*{Supplementary Information}

In the Appendix, we present additional numerical results and properties of the mathematical objects used in the main text. We provide proofs of less important quantities that support argument of the main text. These materials are not crucial to the main text; however, they provide additional insight into the mechanism of the algorithm and the meaning and relevance of the functions used in numerical simulations. 

In Appendix \ref{app:numerics}, we give supplementary numerical results: a typical example where the sketch and the rules present a mismatch and a comparaison of matching pursuit against basis pursuit used as decoder.

In Appendix \ref{app:coherence}, we discuss the coherence we mentioned in the Theoretical Analysis. In particular, we give a formal expression for the coherence of structured matrices that were used in the numerical simulations of the Numerical Results. 

In Appendix \ref{app:theory_check}, we test the validity of the bound provided by Theorem 1 from the main text (Theoretical Analysis) with simulations.

In Appendix \ref{app:threshold}, we give the proofs of the results on thresholding stated in the Theoretical Analysis, namely Proposition 1 and Theorem 2 of the main text, together with the two lemmas they rely on.

\section{Supplementary numerical results}
\label{app:numerics}

\begin{figure}[ht]
	\centering
	\includegraphics[width=0.7\linewidth]{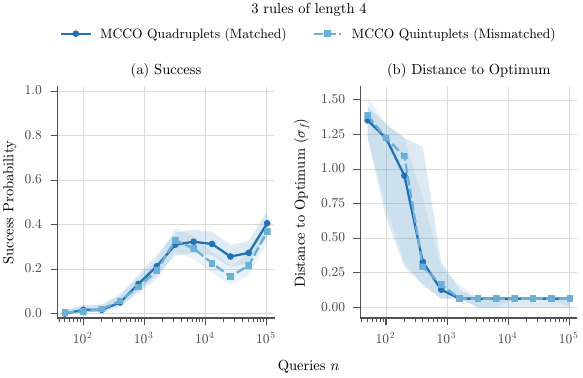}
	\captionsetup{justification=raggedright}
	\caption{Sketch matched or mismatched to the rule length, on one instance with $3$ local rules of length $4$. (a) Success probability and (b) distance to the optimum $(f(x^*) - f(\hat{x}))/\sigma_f$ as a function of the number of queries $n$, over $300$ runs, for the quadruplet sketch (matched) and the quintuplet sketch (mismatched). Bands: $95\%$ Wilson intervals in (a), interquartile range around the median in (b).}
	\label{fig:e3_mismatch}
\end{figure}

\begin{figure}[ht]
	\centering
	\includegraphics[width=0.8\linewidth]{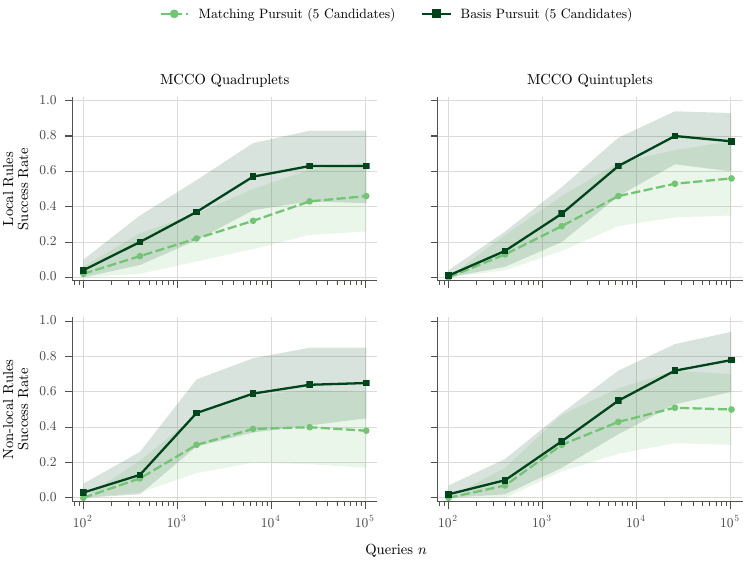}
	\captionsetup{justification=raggedright}
	\caption{Matching pursuit and basis pursuit as decoder, both returning $5$ candidates, for the quadruplet and quintuplet sketches. Success rate as a function of the number of queries $n$ on $20$ instances per family with $|\mathcal{R}| = 5$ rules, $5$ runs per instance, local rules (top row) and non-local rules (bottom row). Shaded bands are $95\%$ intervals.}
	\label{fig:s1_decoders}
\end{figure}

This section gathers two numerical results that complete the Numerical Results of the main text. Supplementary Figure~\ref{fig:e3_mismatch} tests whether the sketch has to match the rules of the objective, and Supplementary Figure~\ref{fig:s1_decoders} compares matching pursuit with basis pursuit as decoding procedure. Unless stated otherwise, the instances, the sampling protocol, the threshold and the estimate are those of the main text.

Supplementary Figure~\ref{fig:e3_mismatch} isolates the effect of a mismatch between the sketch and the rules. The instance is typical and drawn from the local family with $3$ rules, all of length $4$, so that the Quadruplets sketch matches the length of the rules while the Quintuplets sketch does not. Over $300$ runs, which means $300$ different sample choices, the two sketches reach comparable success probabilities at most budgets, e.g. $0.31$ and $0.33$ at $n = 3200$ and $0.41$ and $0.37$ at $n = 102400$, and their median distance to the optimum decreases in the same way. The mismatched sketch therefore performs as well as the matched one: we conclude that MCCO does not need to know the exact length of the rules.

We compare matching pursuit with basis pursuit as decoding procedure, on $20$ instances per family with $|\mathcal{R}| = 5$ and $5$ runs per instance, both decoders returning $5$ candidates (Supplementary Figure~\ref{fig:s1_decoders}). The success rates agree within the confidence intervals up to $n = 1600$. At larger budgets, basis pursuit reaches higher success rates ($0.63$--$0.78$ against $0.38$--$0.56$ at $n = 102400$), at a decoding time of $33$--$41$\,s against $0.01$--$0.02$\,s for matching pursuit, i.e. more than three orders of magnitude.

These results complete the picture of the main text. MCCO does not require the sketch to match the rules of the objective. Matching pursuit performs as well as basis pursuit when queries are few, which is the regime of interest, at a decoding cost more than three orders of magnitude lower.

\newpage
\section{Coherence of sketch matrices}
\label{app:coherence}

\begin{theorem}
	The coherence of a structured matrix using nearest-neighbor patterns of size $k$ acting on a space of $N$ bits is given by
	$$
	\mu = \frac{N-k}{N-k+1}
	$$
\end{theorem}
\begin{proof}
	The $k$-uplet can be positioned starting from $N-(k-1)$ distinct positions (e.g. $N-1$ for pairs, $N-2$ for triplets, etc.).
	
	So any column of the structured matrix has $N-k+1$ entries and the associated normalization factor is $\frac1{\sqrt{N-k+1}}$.
	
	Then, for each of the $N-k+1$ starting positions, the first column matches with the second column except for the last position, i.e. when we look at each bit individually.
	
	The number of matches that contribute to the scalar product is then $N-k$.
\end{proof}

\begin{theorem}
	The coherence of a structured matrix using all possible patterns of size $k$ acting on a space of $N$ bits is given by
	$$
	\mu = \frac{\binom{N-1}{k}}{\binom{N}{k}}
	$$
\end{theorem}
\begin{proof}
	For full matrices, we can place the measured $k$-uplet at $\binom{N}{k}$ starting positions.
	
	So any column of the structured matrix has $\binom{N}{k}$ entries and the normalization factor is $\frac1{\sqrt{\binom{N}{k}}}$.
	
	Then, for each of the $\binom{N}{k}$ starting positions, the first column matches each time with the second column except when one of the measured bits is in the last position, i.e. when we look at each bit individually.
	
	The number of cases where this does not happen is exactly $\binom{N-1}{k}$.
\end{proof}

For the random Gaussian sketch used in the Numerical Results ($512$ rows and $2^{20}$ columns), the coherence cannot be computed exactly, as it is a maximum over about $5 \times 10^{11}$ pairs of columns. We computed it on random subsets of $2^{15}$ to $2^{16}$ columns and extrapolated it to $2^{20}$ columns with the extreme-value scaling of the maximum of Gaussian variables, which gives a coherence of about $0.31$.

\newpage
\section{Numerical verification of the theory}
\label{app:theory_check}

In this section, we test Theorem 1 of the main text on instances of the Numerical Results. The theorem bounds the probability that Problem II fails, i.e. that the maximizer of the empirical surrogate $\tilde{F}_{S_n}$ differs from $x^*$, by Eq.~(6), and gives the corresponding sample size, Eq.~(7). We compute these quantities exactly, with the threshold $t$ set to the $85$th percentile of $f$, and compare them with the empirical failure rate of Problem II.

\begin{figure}[ht]
	\centering
	\includegraphics[width=1\linewidth]{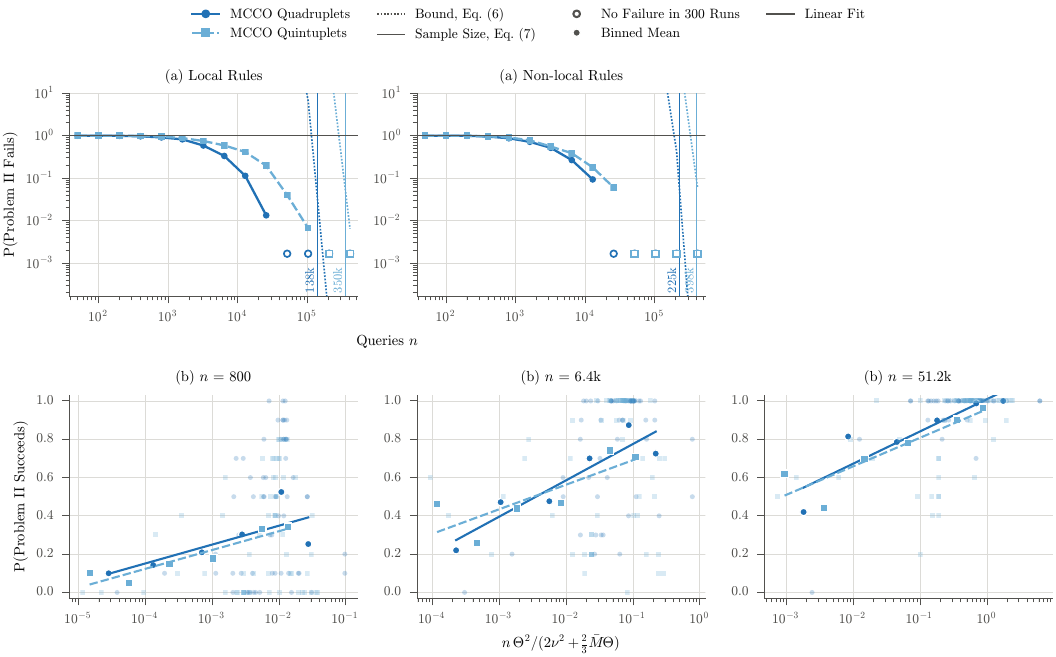}
	\captionsetup{justification=raggedright}
	\caption{Comparison with Theorem 1 of the main text, with the threshold $t$ set to the $85$th percentile of $f$. \\
    (a) Empirical probability that Problem II fails, $\argmax \tilde{F}_{S_n} \neq x^*$, as a function of $n$, over $300$ runs, on one instance per family with a unique maximizer. The dotted lines show the theoretical bound of Eq.~(6); the vertical lines represent the sample size given vby Eq.~(7) for $\delta = 0.1$; the open markers stands for no failure in $300$ runs (zero can't be drawn on the scale). \\
    (b) Per-instance success probability of Problem II over $10$ sample choices, on the instances of Figure 3 of the main text, as a function of the exponent $n\Theta^2/(2\nu^2 + \frac{2}{3}\bar{M}\Theta)$ of Eq.~(6), at $n = 800$, $6400$ and $51200$. The light points show single instances succeeds rate; the filled markers are mean values; finally the trend curvers are computed from linear fits of the means.}
	\label{fig:e5_theory}
\end{figure}

Supplementary Figure~\ref{fig:e5_theory}(a) shows, on two instances with a unique maximizer for which the bound is informative within our budgets, that the empirical probability that Problem II fails decays exponentially with $n$ and vanishes in $300$ runs from $n = 25600$--$204800$. The bound is conservative by a factor of $2$ to $9$, as expected from the union bound over the $2^N$ competitors, but it orders the sketches correctly: Quadruplets require fewer samples than Quintuplets on both instances, in the simulations as in the bound. Supplementary Figure~\ref{fig:e5_theory}(b) checks the bound across instances: the success of Problem II increases with the exponent of Eq.~(6); averaged over instances with similar exponents (filled markers), the success probability grows by $0.1$ to $0.2$ for every tenfold increase of the exponent, at each of the three sample sizes. A trend curve computed by linear fit show this tendency.

The bound is therefore conservative but informative: it predicts the exponential decay of the failure probability with the number of queries and the ordering of the sketches, and its exponent correlates with the difficulty of the instances.

\newpage
\section{Effect of thresholding}
\label{app:threshold}

We use the notation of the main text: $x^*$ is the unique maximizer of $f$, $x_2$ its second maximum, $T_t$ the thresholding operator, and $\Delta_{x,t}$, $\Theta_{x,t}$, $\nu_{x,t}^2$, $\Theta_{\min,t}$, $\nu_t^2$ and $M_t$ are the quantities of Definition 2.

\begin{proof}[Proof of Proposition 1]
	Since $f \geq 0$, thresholding only replaces values by zero, so $T_t f(x) \leq f(x)$ for all $x$.
	As $t \leq f(x_2)$, both the maximum value and the second largest value survive: $T_t f(x^*) = f(x^*)$ and $T_t f(x_2) = f(x_2)$. \\
	For any $x \neq x^*$, $T_t f(x) \leq f(x) \leq f(x_2)$. Hence $\max_{x \neq x^*} T_t f(x) = f(x_2) < T_t f(x^*)$, and the gap of $T_t f$ is $f(x^*) - f(x_2) = \Delta$.
\end{proof}

\begin{lemma}[Monotonicity of the second moment]
	\label{lem:threshold_moment}
	For every $x \neq x^*$, $\nu_{x,t}^2$ decreases as $t$ increases, and hence so does $\nu_t^2 = \max_{x \neq x^*} \nu_{x,t}^2$.
\end{lemma}
\begin{proof}
	We use the fact that $T_t f(s) = \mathbf{1}_{\{f(s) \geq t\}}\, f(s)$ where $\mathbf{1}$ is the indicator function. Taking the difference, we have $\Delta_{x,t}(s)^2 = \mathbf{1}_{\{f(s) \geq t\}}\, \Delta_x(s)^2$, and the indicator can only switch from $1$ to $0$ as $t$ increases.
\end{proof}

\begin{definition}[Discarded part, discarded surrogate and discarded gap]
	\label{def:discarded}
	For a threshold $t \geq 0$ and $x \neq x^*$, the part of $f$ discarded by thresholding, its surrogate, and its surrogate gap are
	\begin{equation*}
		(T_t f)^\complement = f - T_t f, \qquad F^\complement_t(x) = 2^{-N}\left(G\,(T_t f)^\complement\right)_x, \qquad \Theta^\complement_{x,t} = F^\complement_t(x^*) - F^\complement_t(x).
	\end{equation*}
\end{definition}

\begin{lemma}[Surrogate gap under thresholding]
	\label{lem:threshold_gap}
	For every $x \neq x^*$,
	\begin{equation}
		\label{eq:thresh_identity}
		\Theta_{x,t} = \Theta_x - \Theta^\complement_{x,t}.
	\end{equation}
	In particular, if $x^*$ is a minimizer of $F^\complement_t$, then
	\begin{enumerate}[label=(\roman*)]
		\item $\Theta_{x,t} \geq \Theta_x$ for every $x \neq x^*$, so that $\Theta_{\min,t} \geq \Theta_{\min}$;
		\item if $G$ preserves the maximum of $f$, it also preserves the maximum of $T_t f$.
	\end{enumerate}
\end{lemma}
\begin{proof}
	The surrogate is linear in the original function it is computed from, so $f = T_t f + (T_t f)^\complement$ gives $F = F_t + F^\complement_t$. 
	Taking the difference between $x^*$ and $x$ gives the identity (\ref{eq:thresh_identity}).\\ 
	If $x^*$ minimizes $F^\complement_t$, then all expected gaps are below zero, i.e. $\Theta^\complement_{x,t} \leq 0$ for all $x \neq x^*$, 
	which gives (i). \\
	If moreover $\Theta_x > 0$ for all $x \neq x^*$, which is the case when $G$ preserves the maximum of $f$, then $\Theta_{x,t} > 0$, which gives (ii).
\end{proof}

The surrogate gap of $f$ splits into the gap of what is kept and the gap of what is discarded, 
and thresholding widens it when the discarded part has $x^*$ as its least favored point. 
This happens because the objective function is assumed to be compressible, through a selection effect: points that resemble $x^*$ share its rewarding patterns, 
so they tend to take large values and are kept, while the discarded set, on the other hand, is depleted of such points.

Conversely, only the kept points contribute to $F_t$, so that $\Theta_{\min,t}$ is at most $p_1(t)\, f(x^*) \max_{s,x}\left|g^{x^*}(s) - g^x(s)\right|$, 
where $p_1(t)$ is the fraction of points $x$ with $f(x) \geq t$. 
As a consequence, the surrogate gap vanishes as the kept fraction decreases. 
In practice, $t$ is chosen as a quantile of the sampled values, which fixes $p_1(t)$ and hence this bound.

\begin{proof}[Proof of Theorem 2]
	By Lemma \ref{lem:threshold_gap}, $G$ preserves the maximum of $T_t f$. This is the condition for Corollary 1 of the main text to apply, and, by Lemma \ref{lem:threshold_gap} again, thresholding increases the expected surrogate gap, i.e. $\Theta_{\min,t} \geq \Theta_{\min}$. \\
	In addition, Lemma \ref{lem:threshold_moment} shows that the second moment decreases with thresholding, i.e. $\nu_t^2 \leq \nu^2$. \\
	Since $\|T_t f\|_\infty = \|f\|_\infty$, the bound on $M$ in the proof of Theorem 1 gives $M \leq \bar{M}$ 
	and $M_t \leq \bar{M}$, so Theorem 1 holds with $\bar{M}$ for both $f$ and $T_t f$. \\
	Finally, the map $(u,v) \mapsto u^2/(2v+\frac{2}{3}\bar{M}u)$ is increasing in $u > 0$ and decreasing in $v$, which proves the inequality of Theorem 2.
\end{proof}

\end{document}

%% file: fig/e4_cost.tex
\begin{tabular}{llrrrrr}
\hline
Method & Family & sampling & sketching & decoding & candidates & total \\
 & & \multicolumn{5}{c}{median wall-clock (s) at $n = 6400$} \\
\hline
MCCO Quadruplets & Local & 0.151 & 0.0144 & 0.0113 & 0.000373 & 0.178 \\
MCCO Quadruplets & Non-Local & 0.151 & 0.0139 & 0.011 & 0.000363 & 0.178 \\
MCCO Quintuplets & Local & 0.151 & 0.0212 & 0.0181 & 0.000371 & 0.192 \\
MCCO Quintuplets & Non-Local & 0.151 & 0.0205 & 0.0178 & 0.000365 & 0.191 \\
MCCO Random & Local & 0.153 & 0.0845 & 11.5 & 0.000432 & 11.9 \\
MCCO Random & Non-Local & 0.151 & 0.077 & 11 & 0.000402 & 11.4 \\
Digital Annealing & Local & -- & -- & -- & -- & 0.00705 \\
Digital Annealing & Non-Local & -- & -- & -- & -- & 0.00706 \\
\hline
\end{tabular}